\documentclass[11pt]{article}
\usepackage{amssymb,amsfonts,amsmath,latexsym,epsf,tikz,url}

\newtheorem{theorem}{Theorem}[section]

\newtheorem{corollary}[theorem]{Corollary}
\newtheorem{lemma}[theorem]{Lemma}
\newtheorem{define}[theorem]{Definition} 

\newtheorem{remark}[theorem]{Remark}

\newcommand{\proof}{\noindent{\bf Proof.\ }}
\newcommand{\qed}{\hfill $\square$\medskip}

\begin{document}

\title{Domination polynomial of clique cover product of graphs}

\author{Somayeh Jahari$^{}$\footnote{Corresponding author}
\and
Saeid Alikhani 
}

\date{}

\maketitle

\begin{center}
 Department of Mathematics, Yazd University, 89195-741, Yazd, Iran
{\tt s.jahari@gmail.com ~~~~ alikhani@yazd.ac.ir}
\end{center}


\begin{abstract}
Let $G$ be a simple graph of order $n$.
The domination polynomial of $G$ is the polynomial
$D(G, x)=\sum_{i=1}^n d(G,i) x^i$,
where $d(G,i)$ is the number of dominating sets of $G$ of size $i$.  
 For two graphs $G$ and $H$, let $\mathcal{C} = \{C_1,C_2, \cdots, C_k\}$ be  a clique cover of $G$ and $U\subseteq V(H)$. We consider clique cover product which denoted by $G^\mathcal{C} \star H^U$ and obtained   from $G$ as follows:  for each clique $C_i \in \mathcal{C}$, add a  copy of the graph $H$ and join every vertex of $C_i$ to every vertex of $U$. 
 We prove that   the domination polynomial of clique cover product $G^\mathcal{C} \star H^{V(H)}$ or simply $G^\mathcal{C} \star H$ is 
 \[ 
	D(G^\mathcal{C} \star H,x)=\prod_{i=1}^k\Big [\big((1+x)^{n_i}-1\big)(1+x)^{|V(H)|}+D(H,x)\Big],
\]
where each clique $C_i \in \mathcal{C}$ has  $n_i$ vertices.  
 As results, we study the  $\mathcal{D}$-equivalence classes of some families of  graphs. Also  we completely describe 
the  $\mathcal{D}$-equivalence classes of friendship graphs constructed by coalescence $n$ 
copies of the cycle graph 
 of length three with a common vertex. 
\end{abstract}


\medskip
\noindent{\bf AMS Subj. Class.:} 05C60, 05C69

\section{Introduction}

All graphs in this paper are simple of finite orders, i.e., graphs are undirected with no loops or
parallel edges and with finite number of vertices. Graph polynomials
are a well-developed area useful for analyzing properties of graphs.  
Li and Gutman \cite{44} introduced a general graph polynomial.
Let $f$ be a complex-valued function defined on the set of graphs $G$ such that $G_1\sim G_2$ 
implies $f(G_1) = f(G_2)$. Let $G$ be a graph on $n$ vertices and $S(G)$ be the set of all subgraphs of $G$. Define 
$S_k(G)=\{H: H \in S(G)  ~\textsl{and}  ~ |V (H)| = k\}$, $p(G, k) = \sum_{H\in S_k(G)}f(H)$.
Then, the general graph polynomial of $G$ is defined as
$P(G,x) =\sum_{k=0}^np(G,k)x^{k}$.
Non-isomorphic graphs may have the same graph polynomial. 
Two graphs
$G$ and $H$ are said to be $P$-equivalent, written as $G \sim^P H$, if $P(G, x) = P(H, x)$. The $P$-equivalence class of $G$ is $[G] = \{H: H\sim^P G\}$. A graph $G$ is said to
be $P$-unique if $[G] = \{G\}$. There are two interesting problems
on the equivalence classes:

(i) Determine the $P$-equivalence classes for some families of graphs.

(ii) Which graphs are $P$-unique?  

These problems have been widely studied for chromatic polynomial (see for example, \cite{Dong}). Domination polynomial of graph $G$ is the  generating function for the number of dominating sets of  $G$, i.e., $D(G,x)=\sum_{ i=1}^{|V(G)|} d(G,i) x^{i}$ (see \cite{euro,saeid1}).
Calculating the domination polynomial of a graph $G$ is difficult in general, as the smallest power of a non-zero term is the domination number $\gamma (G)$ of the graph, and determining whether $\gamma (G) \leq k$ is known to be NP-complete \cite{garey}. But for certain classes of graphs, we can find a closed form expression for the domination polynomial.  The  equivalence classes of the domination polynomial  is called ${\cal D}$-equivalence class.  
 It is known
that cycles \cite{euro} and cubic graphs of order $10$ \cite{cubic} (particularly, the Petersen
graph) are ${\cal D}$-unique, while if $n\equiv 0  (mod\, 3)$, the paths of order $n$ are not \cite{euro}. In \cite{complete}, authors  
have given a necessary and sufficient condition for $\mathcal{D}$-uniqueness of the complete $r$-partite graphs. 
Their results in the bipartite case, settles in the affirmative a conjecture in \cite{ghodrat}. In \cite{barbell} the $\mathcal{D}$-equivalence class of barbell graph (and its generalization) described and showed that there are many families of connected graphs in  $\mathcal{D}$-equivalence class of $nK_r$, where $nK_r$ is disjoint union of $n$ complete graph $K_r$. 
The join $G+H$ of two graph $G$ and $H$ with disjoint vertex sets $V(G)$ and $V(H)$ and
edge sets $E(G)$ and $E(H)$ is the graph union $G\cup H$ together with all the edges joining $V(G)$ and
$V(H)$. For two graphs $G = (V,E)$ and $H=(W,F)$, the corona $G\circ H$ is the graph arising from the
disjoint union of $G$ with $| V |$ copies of $H$, by adding edges between
the $i$th vertex of $G$ and all vertices of $i$th copy of $H$ \cite{Fruc}. It is easy to see that the corona operation of two graphs does not have the commutative property.  
 The induced subgraph $\langle U\rangle$ is a graph with the vertex set $U$ and the edge set consist of edges in $G$ which connect vertices in $U$, if $U \subseteq V (G)$.   A clique is a subset of vertices of an undirected graph such that its induced subgraph is complete; that is, every two distinct vertices in the clique are adjacent.  
 Zhu in \cite{bx} defined an operation of graphs called the clique cover product.  
Given two graphs $G$ and $H$, assume that $\mathcal{C} = \{C_1,C_2, \cdots ,C_k\}$ is a clique cover of $G$ and $U$ is a subset of $V(H)$. Construct a new graph from $G$, as follows: for each clique $C_i \in \mathcal{C}$, add a
copy of the graph $H$ and join every vertex of $C_i$ to every vertex of $U$. Let $G^\mathcal{C} \star H^U$ denote
the new graph. In fact, the clique cover product of graphs is a common generalization of
some known operations of graphs. For instance: If each clique $C_i$ of the clique cover $\mathcal{C}$ is
a vertex, then $G^{V(G)} \star H^{V(H)}$ is the corona of $G$ and $H$. If we take $H = 2K_1$ and $U = V(2K_1)$, then $G^\mathcal{C} \star H^U$ is the graph $\mathcal{C}\{G\}$ obtained by Stevanović \cite{stev} using the clique cover construction.
It has proved that   independence polynomial of $G^\mathcal{C} \star H^U$  is 
\[
I(G^\mathcal{C} \star H^U)= (I(H,x))^k I\Big(G,\frac{xI(H-U,x)}{I(H,x)}\Big),
\]
and showed that $I(G^\mathcal{C} \star H^U)$ is  unimodal \cite{bx}.
We prove that   the domination polynomial of clique cover product $G^\mathcal{C} \star H^{V(H)}$ or simply $G^\mathcal{C} \star H$ is 
\[ 
D(G^\mathcal{C} \star H,x)=\prod_{i=1}^k D(K_{n_i}+H,x),
	\]
where $K_{n_i}=\langle C_i \rangle$.  

 In the next section, we consider graphs which constructed from the path
 $P_n$ by the clique cover construction and study their domination polynomials and ${\cal D}$-equivalence classes.  In Section 3, We extend the results of Section 2 and obtain domination polynomial of clique cover product of graphs, and  
as some consequences, we determine graphs in the class of some specific $k$-trees.  We completely describe graphs in the class $[F_n]$, where $F_n=K_1+nK_2$ is friendship graph, in Section 4.
 \section{ $\mathcal{D}$-equivalence class of a family of graphs}
 
  In this section, we investigate the $\mathcal{D}$-equivalence classes of a family  of graphs. We state the following definition  (see \cite{Lev,stev}).
  
  \begin{define}
A clique cover of a graph $G$ is a spanning subgraph of $G$, each
component of which is a clique. If $\mathcal{C}=\{C_1,C_2,\dots ,C_q\}$ is a clique cover of $G$,
construct a new graph $H$ from $G$, which is denoted by $H=\mathcal{C}\{G\}$,
as follows: for each clique $C\in \mathcal{C}$, add two new non-adjacent vertices and join them to all the vertices of $C$. Note that all old edge of
$G$ are kept in $\mathcal{C}\{G\}$. 
\end{define} 

\begin{figure}[!ht]
\hspace{2.6cm}
\includegraphics[width=9.3cm,height=1.6cm]{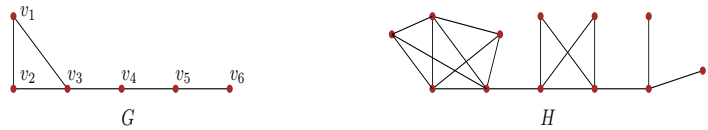}
\caption{ \label{examcliq} Graphs $G$ and $H=\mathcal{C}\{G\}$, respectively.}
\end{figure}

In Figure \ref{examcliq}, the set  $\mathcal{C} =\{\{v_1, v_2, v_3\},\{v_4, v_5\}, \{v_6\}\}$ is a clique cover of $G$ that has a clique consisting of one vertex.

\begin{figure}[!ht]
\hspace{1.9cm}
\includegraphics[width=11cm,height=4.1cm]{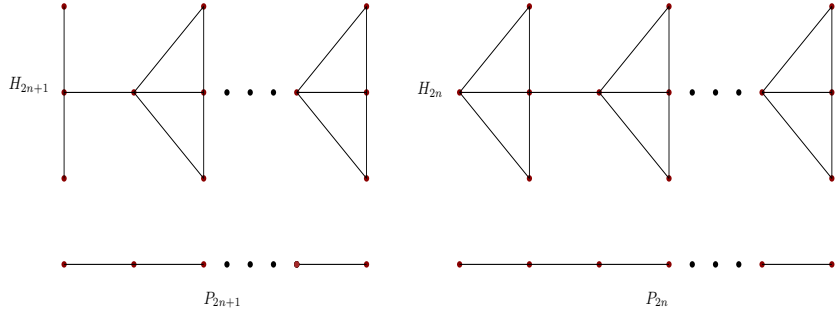}
\caption{ \label{figure1} Graphs $H_{2n+1}$ and $H_{2n}$, respectively. }
\end{figure}

Now,  we  consider graphs of the form $H_n=\mathcal{C}\{P_n\}$, which constructed from the path $P_n$ by
the clique cover construction. Note that in $H_n=\mathcal{C}\{P_n\}$ (Figure \ref{figure1}), for even $n$,  $\mathcal{C} =\{\{1,2\}, \{3,4\}, ..., \{n-1,n\}\}$, and
for odd $n$, $\mathcal{C} = \{\{1\},\{2,3\}, ..., \{n-3,n-2\},\{n-1,n\}\}$.  By $H_0$ we mean the null graph. We shall study the $\mathcal{D}$-equivalence class of $H_n$.
To compute the domination polynomial of $H_n$,  we  need some preliminaries and well known results.

An approach to computing the domination polynomial of
a graph is in term of those of its subgraphs. For instance, one can deduce that 
\[D(G_1\cup G_2,x)= D(G_1,x)D(G_2,x),\]
and the following theorem which gives the domination polynomial of join of two graphs.

\begin{theorem}\label{join}{\rm\cite{euro}}
	Let $G_1$ and $G_2$ be  graphs of orders $n_1$ and $n_2$,
	respectively. Then
	\[
	D(G_1+ G_2,x)=\Big((1+x)^{n_1}-1\Big)\Big((1+x)^{n_2}-1\Big)+D(G_1,x)+D(G_2,x).
	\]
\end{theorem}

The following theorem which is for computation of domination
polynomial of corona products of two graphs.

\begin{theorem}\label{thdcorona}{\rm \cite{Oper,Kot}}
Let $G = (V,E)$ and $H=(W,F)$ be nonempty graphs of order $n$ and $m$, respectively. Then
\begin{eqnarray*}
D(G\circ H,x) = (x(1 + x)^m + D(H, x))^n.
\end{eqnarray*}
\end{theorem}

 We also denote by $G-e$ the subgraph of $G$, obtained by deleting an edge $e$ of $E(G)$. 
  An irrelevant edge is an edge $e\in E(G)$, such that $D(G, x) = D(G-e, x)$, and a vertex $v \in V(G)$ is domination-covered,
if every dominating set of $G-v$ includes at least one vertex adjacent to $v$ in $G$ \cite{Kot}. The following theorem gives a necessary and sufficient condition for a vertex to be a domination-covered vertex.  

\begin{theorem}{\rm\cite{Kot}}\label{dcov}
Let $G = (V,E)$ be a graph. A vertex  $ v\in V$ is domination-covered if and only if there is a vertex  $u\in N[v]$ such that $N[u] \subseteq N[v]$.
\end{theorem}

Using Theorem \ref{dcov} we are able to determine an irrelevant edge. 
\begin{theorem}{\rm\cite{Kot}}\label{irre}
Let $G = (V,E)$ be a graph. An edge $e = \{u, v\} \in E$ is an irrelevant edge
in $G$, if and only if $u$ and $v$ are domination-covered in $G-e$.
\end{theorem}

 Now, we are ready to use Theorem \ref{irre} to obtain  the domination polynomials of $H_n$:
\begin{theorem}\label{thdhn}
 Let   $H_n$ be the  graphs in the Figure \ref{figure1}.
 \begin{enumerate}
 \item[(i)] For every $n \in \mathbb{N}$, $D(H_{2n},x)=(x^4+4x^3+6x^2+2x)^n.$
 \item[(ii)]  For every $n \in \mathbb{N}$, $D(H_{2n+1},x)=(x^3+3x^2+x)(x^4+4x^3+6x^2+2x)^n.$
 \end{enumerate}
\end{theorem}
\proof
\begin{enumerate}
 \item[(i)] Let $G=K_1+ P_3$ be a graph of order $4$
  and  $e_1,\dots,e_n$ be the edges with end-vertices of degree $4$, whose   connect each  two $G$ in $H_{2n}$.
By Theorem \ref{dcov} two end-vertices of every edge $e_i$ are  domination-covered in $H_{2n}$, and so by Theorem \ref{irre}  every edge $e_i$ is an  irrelevant edge of $H_{2n}$.  Since 
 $D(G,x)=x^4+4x^3+6x^2+2x$,  using  induction we have
 $$D(H_{2n},x)=(x^4+4x^3+6x^2+2x)^n. $$
 \item[(ii)] Let $e$ be an edge  joining $H_{2n}$ and $P_3$ in $H_{2n+1}$. By Theorem \ref{dcov} two end-vertices of edge $e$ are  domination-covered in $H_{2n+1}$. So, by Theorem \ref{irre}  the edge $e$ is an  irrelevant edge of $H_{2n+1}$.  Therefore 
 $D(H_{2n+1},x)=  D(P_3\cup H_{2n},x)$ and by Part $(i)$ we have the result.\quad\qed
\end{enumerate}

 The following corollary is an immediate consequence of Theorem \ref{thdhn}.
 \begin{corollary}
$(i)$ For each natural number $n$, the  graph $H_{2n}$ is not $\mathcal{D}$-unique.\\
 $(ii)$  For each natural number $n$, the  graph $H_{2n+1}$ is not $\mathcal{D}$-unique.
\end{corollary}
{\bf Proof.}
$(i)$ Let $G=K_1+ P_3$ be the graph of order four. Since $D(G,x)=x^4+4x^3+6x^2+2x$, so for  $H=n G$, we have 
  $D(H,x)=(x^4+4x^3+6x^2+2x)^n=D(H_{2n},x)$.\\
$(ii)$ We know  $D(H_{2n+1},x)=  D(P_3\cup H_{2n},x)$. If consider  the graph $H$   in the proof of Part $(i)$, then  $P_3\cup H \in [H_{2n+1}]$. \quad\qed

 The following lemma present  many  graphs in the class $[H_{2n}]$:

\begin{lemma}
Let $G$ be a graph of order $n$. The graphs of the form $G\circ P_3$ and $H_{2n}$ have the same domination polynomial. 
\end{lemma}
{\bf Proof.} 
 By theorem \ref{thdcorona} we can deduce that for each arbitrary graph $G$,
\begin{eqnarray*}
D(G\circ P_3,x)&=&\Big(x(1 + x)^3 + D(P_3, x)\Big)^{|V(G)|}\\
&=& (x(1 + x)^3 + x^3+3x^2+x )^{n}\\
&=&\Big(x(x^3+4x^2+6x+2)\Big)^{n}\\
&=&D(H_{2n},x). \qquad\qquad \Box
\end{eqnarray*}

\begin{figure}[!ht]
\hspace{2.7cm}
 \includegraphics[width=9cm,height=4.cm]{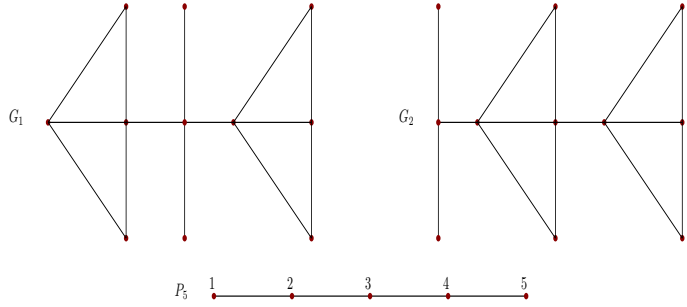}
\caption{ \label{nisoph2} The graphs $G_1=\mathcal{C}_1\{P_{5}\}$ and $G_2=\mathcal{C}_2\{P_{5}\}$, respectively.  }
 \end{figure}

One can see that graphs constructed by clique cover of a graph $G$, may are not isomorphic. As instance,  two  non-isomorphic graphs $G_1$ and $G_2$, dipicted in Figure \ref{nisoph2}, are obtained by different clique covers of $P_5$, namely $\mathcal{C}_1=\{\{1,2\}, \{3\}, \{4,5\}\}$ and $\mathcal{C}_2=\{\{1\}, \{2,3\}, \{4,5\}\}$.
We  notice that the domination polynomials of this two graphs are the same, i.e., 
\[ 
D(G_1,x)=D(G_2,x) = (x^4+4x^3+6x^2+2x)^2(x^3+3x^2+x).
\]

Finally,   we  present some other families of graphs whose are in the $\mathcal{D}$-equivalence classes of $H_n$ graphs in Figure \ref{figure1}. 

\begin{theorem}\label{dunihn}
$(i)$ Let $e_1,\dots, e_n$ be the edges with end-vertices of degree four in the graph $H_{2n}$. For every $1\leq i \leq n$, the disconnected graphs obtained from deletion of any number of the edges $e_i$  and the graph $H_{2n}$ have the same domination polynomial. As well as all graphs obtained by adding each number of the edges between every two vertices of degree four, and adding the edges between two vertices  of degree three and four in the graph $H_{2n}$.

$(ii)$  Let $\mathcal{C}$ be a clique cover of $P_{2n+1}$ that has a clique consisting of one vertex. 
 All non-isomorphic graphs $\mathcal{C}\{P_{2n+1}\}$ and the graph $H_{2n+1}$ have the same domination polynomial. As well as  the union of the graph $P_3$ and every graph in $[H_{2n}]$. 
\end{theorem}
\proof
\begin{enumerate}
 \item[(i)] Similar  to the proof of Theorem \ref{thdhn},  for every $1\leq i\leq n$,  the edge $e_i$ is an  irrelevant edge of $H_{2n}$, that is for every $1\leq i\leq n$, $D(H_{2n}, x) = D(H_{2n}-e_i, x)$, and so we have the first results.  Remains to show that each added edge with the mentioned conditions in this graph is  an irrelevant edge.  This can be achieved by using Theorems \ref{dcov} and \ref{irre}.
 
 \item[(ii)] Similar to the proof of Part $(ii)$ of Theorem \ref{thdhn} we have the result.\quad\qed
\end{enumerate}


\section{Domination polynomial of clique cover product}

  The domination polynomials of binary graph operations, such as, join and  corona has computed \cite{Oper}. Also, recently,  recurrence formulae and properties of the domination polynomials of families of graphs obtained by
  various products, has investigated \cite{complete}. In this section, we generalize the results in Section 2 and consider the clique cover product and  formulate the domination polynomial for the clique cover product $G^\mathcal{C} \star H^{V(H)}$ or simply $G^\mathcal{C} \star H$.  The following theorem gives the domination polynomial of $G^\mathcal{C} \star H^U$.  
  
\begin{theorem}\label{dcli}
 For two graphs $G$ and $H$, let $\mathcal{C} = \{C_1,C_2, \cdots, C_k\}$ be  a clique cover of $G$ and $U\subseteq V(H)$. Then 
\[ 
	D(G^\mathcal{C} \star H^U,x)=\prod_{i=1}^k D(H^*,x),
\]
where $H^*$ is the subgraph of order $|V(H)| + |C_i|$ in $G^\mathcal{C} \star H^U$ obtained by adding a copy of the graph $H$ and joining every vertex of $C_i$ to every vertex of $U$. 
Moreover, 
	\[ 
	D(G^\mathcal{C} \star H,x)=\prod_{i=1}^k \Big[((1+x)^{n_i}-1)(1+x)^{|V(H)|}+D(H,x)\Big],
\]
where $n_i$ is the order of $C_i$. 
\end{theorem} 
\proof 
Since every vertex of $C_i$ is dominated by vertices every vertex in $U$ of $H$ thus by Theorem \ref{dcov} two end-vertices of all edges whose connect each  two $C_i$ in $G^\mathcal{C} \star H^U$  are  domination-covered in new graph, and so by Theorem \ref{irre}  every edge $e_i$ is an  irrelevant edge of $G^\mathcal{C} \star H^U$. Therefore by definition of an  irrelevant edge, and the Principle of Mathematical Induction we have the result. Now, suppose that  $U=V(H)$. Thus  we can deduce that  $ D(G^\mathcal{C} \star H,x)=\prod_{i=1}^k D(K_{n_i}+H,x)$. Note that in the complete graph $K_n$, any nonempty set of vertices is a dominating set, so it follows that $D(K_n, x) = (1+ x)^n -1$. Therefore by Theorem \ref{join}    have the result. \qed

\begin{remark}
	If each clique $C_i$ of the clique cover $\mathcal{C}$ is
a vertex, then $G^{V(G)} \star H$ is the corona of $G$ and $H$. So  the clique cover product of graphs is a  generalization of corona product and hence by Theorem \ref{dcli}  we have 
\begin{eqnarray*}
D(G\circ H,x)&=&D(G^{V(G)} \star H)\\
&=&\prod_{1}^n[((1+x)^{1}-1)(1+x)^{m}+D(H,x)]\\
&=&\Big(x(1+x)^m+D(H,x)\Big)^n,  
\end{eqnarray*}
which is another approach for proof of Theorem \ref{thdcorona}.
\end{remark} 

Here we shall  apply Theorem \ref{dcli} to get the following results on $\mathcal{D}$-equivalence class of some graphs.

\begin{corollary}
	\begin{enumerate}
		\item[(i)] Let $G$ and $H$ be two graphs and $\mathcal{C}_1$ and $\mathcal{C}_2$ be two clique cover of $G$.   
	  	 If $|\mathcal{C}_1|=|\mathcal{C}_2|=k$, and for every $C_i \in \mathcal{C}_1$ there exists $C_j \in \mathcal{C}_2$  of the same order, then  the graphs   $G^{\mathcal{C}_1} \star H$ and $G^{\mathcal{C}_2} \star H$ have the same domination polynomial, i.e.,  
		\[ 
		D(G^{\mathcal{C}_1} \star H,x)=D(G^{\mathcal{C}_2} \star H,x)=\prod_{i=1}^k D(K_{n_i}+H,x),
	\]
where $K_{n_i}=\langle C_i \rangle$.  
		\item[(ii)]
		Given two graphs $G$ and $H$, assume that $\mathcal{C}$ is a clique cover of $G$. The graph $G^{\mathcal{C}} \star H$ and all graphs obtained from deletion of any number of the edges between two cliques in the graph $G^{\mathcal{C}} \star H$ have the same domination polynomial. As well as all graphs obtained by adding each number of the edges between every two cliques of $\mathcal{C}$ in the graph $G^{\mathcal{C}} \star H$.
	\end{enumerate}
\end{corollary}

As an example of application of clique cover for determining the ${\cal D}$-equivalence classes of some graphs, we determine some graphs in the class of $k$-stars. Let us to recall some preliminaries. The class of $k$-trees is a very important subclass of triangulated graphs. Harary and Palmer \rm\cite{harary} first introduced $2$-trees in 1968. Beineke and Pippert \rm\cite{bein} gave the definition of a $k$-tree in 1969.

\begin{define} For a positive integer $k$, a $k$-tree, denoted by $T^k_n$, is defined recursively as follows: The smallest $k$-tree is the $k$-clique $K_k$. If $G$ is a $k$-tree with $n\geq k$ vertices and a new vertex $v$ of degree $k$ is added and joined to the vertices of a $k$-clique in $G$, then the larger graph is a $k$-tree with $n +1$ vertices.
\end{define} 

A $k$-star, $S_{k,n-k}$, has vertex set $\{v_1, \ldots , v_n\}$ where $\langle\{v_1, v_2, \ldots , v_k\}\rangle   \cong K_k$ and $N(v_i) = \{v_1, \ldots, v_k\}$ for $k + 1 \leq i \leq n$. The authors in \cite{jahari} calculated the  domination polynomials for some $k$-tree related graphs, specially $k$-star graph and investigated domination roots of this graph. 

\begin{theorem}\rm\cite{jahari}
	For every $k\in \mathbb{N}$ and $n> k$,
	\[D(S_{k,n-k}, x) = (1+x)^{n-k}((1+x)^k -1) +x^{n-k}. \]
\end{theorem}

Here, we present graphs whose domination polynomials are $\prod \limits_{i=1}^mD(S_{k_i,n_i-k_i}, x)$.

\begin{theorem}
	Let $G$ be a graph and $\mathcal{C}= \{C_1,C_2, \cdots ,C_m : |C_i|=k_i \}$ is a clique cover of $G$. 
	If $H$ is an empty graph, 
	then the graphs $G^\mathcal{C} \star H$ and disjoint union of $m$,  $k_i$- star have the same domination polynomial.
\end{theorem}
\proof
Since $S_{k,n-k}=K_k + S$, where $S$ is an is an empty graph,
 by applying Theorem \ref{dcli}, we  have the result.\qed

The following lemma present many graphs in $[mS_{k,n-k}]$ whose domination polynomials are $D(S_{k,n-k}, x)^m$. 
\begin{lemma}
	\begin{enumerate}
		\item[(i)] Let $G$ be a graph of order $m$, then the graphs $G^{V(G)}\star S_{k-1, n-k+1}=G\circ S_{k-1, n-k+1}$ and $m S_{k, n-k}$ have the same domination poynomial.
		
		\item[(ii)]  Let $G$ be a graph and $\mathcal{C}= \{C_1,C_2, \cdots ,C_m : |C_i|=k \}$ is a clique cover of $G$. 
		If $H$ is an empty graph of order $n-k$, then the graphs $G^\mathcal{C} \star H$ and all graphs obtained from deletion of any number of the edges joining $C_i$ and $C_j$ of $\mathcal{C}$ in the graph $G^\mathcal{C} \star H$ have the same domination polynomial. As well as all graphs obtained by adding each number of the edges between $C_i$ and $C_j$ of $\mathcal{C}$ in the graph $G^{\mathcal{C}} \star H$.
	\end{enumerate}
\end{lemma}


\section{$\mathcal{D}$-equivalence classes of friendship graphs}

\begin{figure}[!ht]
	\hspace{5.3cm}
	\includegraphics[width=4.cm,height=2.5cm]{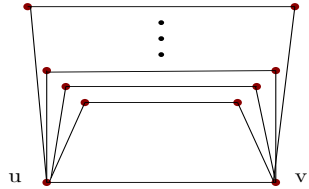}
	\caption{ \label{figure6} The book graph $B_n$.}
\end{figure}

\begin{figure}
	\begin{center}
		\includegraphics[width=5.2in]{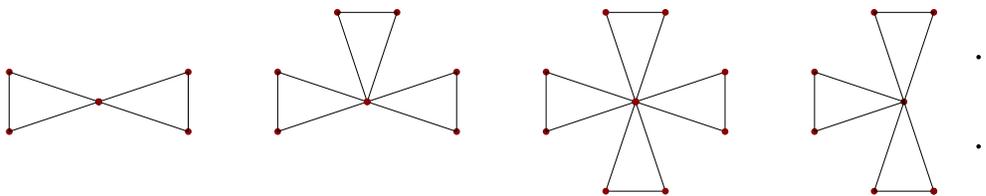}
		\caption{Friendship graphs $F_2, F_3, F_4$ and $F_n$, respectively.}
		\label{Dutch}
	\end{center}
\end{figure}

The friendship (or Dutch-Windmill) graph $F_n$ is a graph that can be constructed by coalescence $n$
copies of the cycle graph $C_3$ of length $3$ with a common vertex. The Friendship Theorem of Paul Erd\"{o}s,
Alfred R\'{e}nyi and Vera T. S\'{o}s \cite{erdos}, states that graphs with the property that every two vertices have
exactly one neighbour in common are exactly the friendship graphs.
Figure \ref{Dutch} shows some examples of friendship graphs. 
The nature and location of domination  roots of friendship graphs have been studied in \cite{FILOMAT} and shown that for every $n\geq 3$, $F_n$ is not $\mathcal{D}$-unique. The authors considered the $n$-book graph $B_n$ which  can be constructed by bonding $n$ copies of the cycle graph $C_4$ along a common edge $\{u, v\}$, see Figure \ref{figure6}. 
The vertex contraction $G/u$ of a graph $G$ by a vertex $u$ is the operation under
which all vertices in $N(u)$ are joined to each other and then $u$ is deleted (see\cite{Wal}). The following result proves that the friendship  graph $F_n$ is not $\mathcal{D}$-unique. 

\begin{theorem}{\rm\cite{FILOMAT}}\label{class} 
	For each natural number $n\geq 3,$ the friendship  graph $F_n$ is not $\mathcal{D}$-unique, as  $F_n$ and $B_n/v$ have the same domination polynomial.
\end{theorem}

  In this section, we describe  
$[F_n]$ completely. Since $F_n=K_1^{V(K_1)}\star (nK_2)=K_1+nK_2$, we have the following theorem.  

\begin{theorem}{\rm \cite{FILOMAT}}
	For every $n\in \mathbb{N}$, $D(F_{n},x) = (2x + x^2)^n +x(1 + x)^{2n}.$
\end{theorem}

We shall  extend Theorem \ref{class} and  present all families of graphs whose are in the $[F_n]$. The following
theorem gives us the domination polynomial of graphs of
the form $H\circ K_1$ which is the first result for domination
polynomial of specific corona of two graphs and we need it to obtain our result.

\begin{theorem}\label{theo}{\rm \cite{euro}}
	 $D(G,x) = x^n(x+2)^n$ if and only if $G=H\circ K_1$
	for some graph $H$ of order $n$.
\end{theorem}

As shown in  \cite{complete}, the following corollary is a consequence of Theorem \ref{join}.

\begin{corollary}\label{corjoin}
	For graphs $G_1, G_2$, and $H$, $D(G_1 + H, x) = D(G_2 + H, x)$ if and only
	if $D(G_1, x) = D(G_2, x)$.
\end{corollary}

The following theorem gives the ${\cal D}$-equivalence classes of $[F_n]$: 
\begin{theorem}
	Let $G$ be a graph of order $n$. Then
	\[
	[F_n]=\{(G\circ K_1)+ K_1 : | G |=n\}.
	\]
\end{theorem} 
\proof
Since $F_n=nK_2+ K_1$, by Corollary \ref{corjoin} to obtain $[F_n]$,   it is suffices  to find $[nK_2]$. Using  Theorem \ref{theo} we have  
\[ 
[nK_2]=\{(G\circ K_1) : | G |=n\}.
\]
So we have the result. \qed

The graph  $B_n/v$ which has found in Theorem \ref{class} is in the form $(K_n\circ  K_1) + K_1$. And 
$F_n$ is one of the graphs in form $(G\circ K_1)+ K_1$ where $G$ is empty graph of order $n$. If $n=1$ then $F_1=K_3$ and by  \cite[Corollary 2]{euro} the  complete graphs $K_n$ are $\mathcal{D}$-unique for every natural number $n$.

\end{document}